\documentclass[reqno]{amsart}
\usepackage{amsmath,amsfonts,epsfig,amssymb,graphics}

\newtheorem{thm}{Theorem}[section]
\newtheorem{corollary}[thm]{Corollary}
\newtheorem{prop}[thm]{Proposition}
\newtheorem{defn}[thm]{Definition}
\newtheorem{remark}[thm]{Remark}
\newtheorem{assumption}[thm]{Assumption}

\begin{document}

\title{Isospectral Finiteness of Hyperbolic Orbisurfaces \ }

\author{Emily B. Dryden}
\address{} 
\email{}
\thanks{Partially supported by NSF grant DMS-0306752. }

\keywords{}

\begin{abstract}
We discuss questions of isospectrality for hyperbolic orbisurfaces, 
examining the relationship between the geometry
of an orbisurface and its Laplace spectrum.  We show that certain hyperbolic
orbisurfaces cannot be isospectral, where the obstructions involve the
number of singular points and genera of our orbisurfaces.  
Using a version of the 
Selberg Trace Formula for hyperbolic orbisurfaces, we show
that the Laplace spectrum determines the length spectrum and the orders 
of the singular points, up to finitely
many possibilities.  Conversely, knowledge of the
length spectrum and the orders of the singular points determines the
Laplace spectrum.  This partial generalization of Huber's theorem is used to prove
that isospectral sets of hyperbolic orbisurfaces have finite cardinality, generalizing 
a result of McKean \cite{McK72} for Riemann surfaces. 
\end{abstract}

\maketitle

\section{Introduction} \label{Sect:S1}
Historically, inverse spectral theory has been concerned with the
relationship between the geometry and the spectrum of compact
Riemannian manifolds, where ``spectrum'' means the eigenvalue spectrum of
the Laplace operator as it acts on smooth functions on a manifold $M$.
Since Milnor's pair of isospectral (same
eigenvalue spectrum) non-isometric flat tori in dimension 16, there
has been much work to characterize the properties of smooth manifolds which
are spectrally determined.  Recently, this study has broadened to 
the category of orbifolds, which are a natural generalization of manifolds 
(see, e.g., \cite{Donheatinvars}, \cite{GGWZ}, \cite{GoRoHodge}, \cite{Stanhopebounds}).  
Instead of being locally modelled on $\mathbb{R}^n$, an orbifold is locally 
modelled on $\mathbb{R}^n$ modulo the action of a discrete group of isometries.   

We will be interested in compact hyperbolic orbisurfaces, which are a natural
generalization of compact hyperbolic Riemann surfaces.  A hyperbolic
orbisurface can be viewed as the quotient of the hyperbolic plane by a finite
group of isometries which is permitted to include elliptic elements.  
These elliptic 
elements give rise to conical singularities in the quotient surface.  
We examine the relationship between the
geometric properties of orbisurfaces and the eigenvalue spectrum of
the Laplace operator as it acts on smooth functions on
the orbisurface.  
For compact Riemann surfaces, Huber's theorem is a powerful tool in studying 
this relationship.  It says that the Laplace spectrum determines the length 
spectrum and vice versa, where the length spectrum is the sequence of lengths 
of all oriented closed geodesics in the surface, arranged in ascending order.  
The Selberg Trace Formula is used to derive this result. 
There has been recent interest in extending Huber's theorem for 
Riemann surfaces to more general settings.  
A trace formula for discrete cocompact groups $\Gamma < PSL(2, \mathbb{C})$ is 
developed by Elstrodt et al. in \cite{ElGrMebook}; they use the formula to 
find a natural generalization of Huber's theorem.  Parnovskii 
\cite{Parnov} states a trace formula in the case of a cocompact discrete 
group of isometries of hyperbolic space of arbitrary dimension; a generalized 
version of Huber's theorem is stated without proof.  In our extension 
of Huber's theorem, we emphasize the geometric meaning of the various 
quantities involved in the Selberg Trace Formula in the orbifold setting.  
Our proof is completely 
different from that given in \cite{ElGrMebook} for dimension three.  
We use our extended Huber's theorem to show that isospectral sets of 
hyperbolic orbisurfaces have finite cardinality, generalizing a result of 
McKean \cite{McK72} for Riemann surfaces.
  
The paper is organized as follows.  
We begin with the necessary definitions and background
concerning orbifolds and the eigenvalue spectrum of the
Laplace operator in the orbifold context.
In section  \ref{sec:obstructions}, 
we move into the study of the relationship between the geometry
of an orbisurface and its Laplace spectrum.  Using Weyl's
asymptotic formula for orbifolds as developed by Farsi \cite{Farsi}, 
we are able to give obstructions to
isospectrality of Riemann orbisurfaces; these obstructions involve the
genera and number of singularities of our orbisurfaces.  We show in section \ref{sec:hubers}
that the Laplace spectrum determines the length spectrum and the orders of the singularities, 
up to finitely
many possibilities.  Conversely, knowledge of the
length spectrum and the orders of the singular points determines the
Laplace spectrum.
In the final section, we use this theorem to show finiteness of
isospectral sets of Riemann orbisurfaces. 

\section*{Acknowledgements}
This work was done while the author was a graduate student at Dartmouth College.  She thanks her advisor, Carolyn Gordon, for her support and for many helpful discussions.  She is also grateful to Peter Buser for his enthusiasm and expertise, and in particular for pointing out that an assumption on the orders of the cone points could be removed from the statements of the theorems in sections \ref{sec:hubers} and \ref{sec:finiteness}.

\section{Preliminaries}\label{sec:prelims}

On orbifolds, we can define coordinate charts which encode information 
about the local group actions at the singularities.  We 
follow Satake \cite{SatakeGB} and Stanhope \cite{Stanhopebounds}.  
Additional references include Chapter 2 of \cite{CoHoKebook}, \cite{Scott83} 
and \cite{Thnotes}.

\begin{defn}
Let $X$ be a Hausdorff space, and let $U$ be an open
set in $X$.  An \emph{orbifold coordinate chart} over $U$ is a triple
($U, \Gamma \backslash \tilde{U}, \pi$) such that:
\begin{enumerate}
\item $\tilde{U}$ is a connected open subset of $\mathbb{R}^n$,
\item $\Gamma$ is a finite group of diffeomorphisms acting on
  $\tilde{U}$ with fixed point set of codimension at least two, and
\item $\pi: \tilde{U} \rightarrow U$ is a continuous map which induces
  a homeomorphism between $\Gamma \backslash \tilde{U}$ and $U$.  We
  require $\pi \circ \gamma = \pi$ for all $\gamma \in \Gamma$.
\end{enumerate}
\end{defn}

Now suppose that $U$ and $U'$ are two open sets in a Hausdorff
space $X$ with $U \subset U'$.  Let ($U, \Gamma \backslash \tilde{U}, \pi$) and
($U', \Gamma' \backslash \tilde{U'}, \pi'$) be charts over $U$ and $U'$,
respectively.  

\begin{defn}\label{defn:injection}
An \emph{injection} $\lambda: (U, \Gamma \backslash \tilde{U}, \pi)
\hookrightarrow (U', \Gamma' \backslash \tilde{U'}, \pi')$ consists of
an open embedding $\lambda:
\tilde{U} \hookrightarrow \tilde{U'}$ such that $\pi = \pi' \circ \lambda$.  
\end{defn}   

\noindent For any $\gamma \in \Gamma$, there exists a unique $\gamma' \in \Gamma'$ for
which $\lambda \circ \gamma = \gamma' \circ \lambda$.  
The correspondence $\gamma \mapsto \gamma'$ defines an
injective homomorphism of groups from $\Gamma$ into $\Gamma'$.

\begin{defn}
A \emph{smooth orbifold} ($X, \mathcal{A}$) consists
of a Hausdorff space $X$
together with an atlas of charts $\mathcal{A}$ satisfying the
following conditions:
\begin{enumerate}
\item For any pair of charts $(U, \Gamma \backslash \tilde{U}, \pi)$ and $(U',
\Gamma' \backslash \tilde{U'}, \pi')$ in $\mathcal{A}$ with $U \subset U'$
there exists an injection $\lambda: (U, \Gamma \backslash \tilde{U}, \pi)
\hookrightarrow (U',\Gamma' \backslash \tilde{U'}, \pi')$.
\item The open sets $U \subset X$ for which there exists a chart $(U,
  \Gamma \backslash \tilde{U}, \pi)$ in $\mathcal{A}$ form a basis of open sets
  in $X$.
\end{enumerate}
\end{defn}

Given an orbifold $(X, \mathcal{A})$, we call the topological space $X$ the
\emph{underlying space} of the orbifold.  Henceforth orbifolds $(X,
\mathcal{A})$ will be denoted simply by $O$.
We now give some examples of orbifolds.

\begin{enumerate}   
\item Let $\Gamma$ be a group acting properly discontinuously on a
  manifold $M$ with fixed point set of codimension at least two.  Then
  the quotient space $O=\Gamma \backslash M$ is an orbifold.  Since
  $O$ can be expressed as a global quotient (that is, as a subset of
  $\mathbb{R}^n$ modulo the action of a discrete group), it is called
  a \emph{good} or \emph{global} orbifold.  If $M$ is a surface, then
  $O$ is an orbisurface.
\item Consider the orbisurface whose underlying space is the sphere
  $S^2$, and which has one singular point.  A neighborhood of this singular
  point is modelled on $\mathbb{Z}_n \backslash \mathbb{R}^2$, where
  $\mathbb{Z}_n$ is the group of rotations of order $n$.  Such a singular point
  is called a \emph{cone point of order $n$}, and this
  orbisurface is known as the $\mathbb{Z}_n$-teardrop.  Unlike the previous
  example, the $\mathbb{Z}_n$-teardrop cannot be 
  expressed as a quotient with respect to the action of a discrete group, 
  and thus is an example of a \emph{bad} orbifold (see \cite{Scott83}). 
\end{enumerate}

We will be interested in orbifolds which have a hyperbolic
structure.  The construction of a Riemannian metric on an orbifold
is as in the manifold case, with the
metric being defined locally via coordinate charts and patched together
using a partition of unity.  In addition, the metric must be invariant
under the local group actions.  A smooth orbifold with a Riemannian
metric is a \emph{Riemannian orbifold}.  An orbisurface with a
hyperbolic metric of constant curvature -1 will be called a
\emph{Riemann orbisurface}.
Every Riemann orbisurface arises as a global quotient of the
hyperbolic plane by a discrete group of isometries (see \cite{Scott83}).  

An orbifold $O$
is said to be \emph{locally orientable} if it has an atlas in which
every coordinate chart $(U, \Gamma \backslash \tilde{U}, \pi)$ is such that
$\Gamma$ is an orientation-preserving group.  If all injections as in
Definition \ref{defn:injection} are
induced by orientation-preserving maps, then $O$ is \emph{orientable}.

Our goal is to study the spectrum of the Laplace operator as it acts on
smooth functions on a compact Riemannian orbifold $O$.  
A map $f: O \rightarrow \mathbb{R}$ is a \emph{smooth function} on $O$
if for every coordinate chart $(U, \Gamma \backslash \tilde{U}, \pi)$ on $O$,
the lifted function $\tilde{f} = f \circ \pi$ is a smooth function on $\tilde{U}$. 
If $O$ is a compact Riemannian orbifold and $f$ is a smooth function
on $O$, then we define the Laplacian $\Delta f$ of $f$ by lifting $f$
to $\tilde{f}$.  We denote the
$\Gamma$-invariant metric on $\tilde{U}$ by $g_{ij}$ and set $\rho =
\sqrt{\emph{det}(g_{ij})}$.  Then we can define the Laplacian locally by 
\begin{displaymath}
\Delta \tilde{f} =
\frac{1}{\rho}\sum_{i,j=1}^{n}\frac{\partial}{\partial \tilde{x^i}}
(g^{ij} \frac{\partial f}{\partial \tilde{x^j}} \rho).
\end{displaymath}

We are really interested in the eigenvalues of the Laplace operator as
it acts on smooth functions.  In analogy with the manifold case,
Chiang \cite{Chiang90} proved the following theorem:

\begin{thm} 
Let $O$ be a compact Riemannian orbifold.  
\begin{enumerate}
\item The set of eigenvalues $\lambda$ in $\Delta f = \lambda f$
  consists of an infinite sequence $0 \leq \lambda_1 <
  \lambda_2 < \lambda_3 < \cdots \uparrow \infty$.  We call this
  sequence the spectrum of the Laplacian on $O$, denoted $Spec(O)$. 
\item Each eigenvalue $\lambda_i$ has finite multiplicity.  
\item There exists an orthonormal basis of $L^2(O)$ composed of smooth
  eigenfunctions $\phi_1, \phi_2, \phi_3, \ldots$, where $\Delta
  \phi_i = \lambda_i \phi_i$.  
\end{enumerate}
\end{thm}

The multiplicity of the $i$th eigenvalue $\lambda_i$ is the dimension
of the space of eigenfunctions with eigenvalue $\lambda_i$.  

\section{Obstructions to Isospectrality}\label{sec:obstructions}

Much information about the relationship between the Laplace spectrum
of an orbifold $O$ and the geometric properties of $O$ can be gleaned by
studying the heat equation:
\begin{displaymath}
\Delta F = - \frac{\partial F}{\partial t},
\end{displaymath}
where $F(x,t)$ is the heat at a point $x \in O$ at time $t$. 
With initial data $f:O \rightarrow \mathbb{R}$, $F(x,0)=f(x)$, a
solution of the heat equation is given by
\begin{displaymath}
F(x) = \int_{O} K(x,y,t)f(y)dy.
\end{displaymath}
Here $K: O \times O \times \mathbb{R}_{+}^{*} \rightarrow \mathbb{R}$
is a $C^{\infty}$ function given by the convergent series
\begin{equation}\label{eqn:heatkernel}
K(x,y,t) = \sum_i e^{- \lambda_i t} \phi_i (x) \phi_i (y).
\end{equation}
The eigenfunctions $\phi_i$ of $\Delta$ are chosen such that they form
an orthonormal basis of $L^2 (O)$, the square-integrable functions on $O$.
We say that $K$ is the \emph{fundamental solution of the heat
  equation} on $O$, or the \emph{heat kernel} on $O$.  The appropriate
  physical interpretation is that $K(x,y,t)$ is the temperature at time $t$
at the point $y$ when a unit of heat (a Dirac delta-function) is
  placed at the point $x$ at time $t=0$.  

By considering the asymptotic behavior of $K$ as $t \rightarrow 0$, we can
recover information about the geometry of $O$.  In this direction,
Farsi showed (see \cite{Farsi}) that Weyl's asymptotic formula can be
extended to orbifolds.  In particular, she proved

\begin{thm}\label{thm:Weylasymp}
Let $O$ be a closed orientable Riemannian orbifold 
with eigenvalue spectrum $0 \leq
\lambda_1 \leq \lambda_2 \leq \lambda_3 \ldots \uparrow \infty$.  Then
for the function $N(\lambda) = \sum_{\lambda_j \leq \lambda} 1$ we
have
\begin{displaymath}
N(\lambda) \sim (\emph{Vol } B_0^n(1))(\emph{Vol }
O)\frac{\lambda^{n/2}}{(2\pi)^n}
\end{displaymath}
as $\lambda \uparrow \infty$.  Here $B_0^n(1)$ denotes the
$n$-dimensional unit ball in Euclidean space, and $n$ is the dimension of $O$.
\end{thm}

\noindent This theorem implies that, in analogy with the manifold case, the
Laplace spectrum determines an orbifold's dimension and volume.  

By looking at the terms of the asymptotic expansion of the trace of
the heat kernel, Gordon et al. \cite{GGWZ} have
given the following obstruction to isospectrality:

\begin{thm}\label{thm:commoncover}
Let $O$ be a Riemannian orbifold with singularities.  If $M$ is a
manifold such that $O$ and $M$ have a common Riemannian cover, then
$M$ and $O$ cannot be isospectral.
\end{thm}

\noindent In particular, this implies that a hyperbolic orbifold with
singularities is never isospectral to a hyperbolic manifold. 

We want to investigate further obstructions to isospectrality; our focus
will be the case of orbisurfaces.  In analogy with the surface case,
we can define the Euler characteristic
and state a Gauss-Bonnet theorem for orbisurfaces (see \cite{Thnotes}).

\begin{defn}
Let $O$ be an orbisurface with $s$ cone points of orders
$m_1,\ldots,m_s$.  Then we define the (orbifold) Euler characteristic
of $O$ to be 
\begin{displaymath}
\chi(O) = \chi (X_0) - \sum_{j=1}^s (1-\frac{1}{m_j}),
\end{displaymath}
where $\chi(X_O)$ is the Euler characteristic of the underlying space of $O$.
\end{defn}

The Gauss-Bonnet theorem for orbisurfaces gives the usual relationship
between topology and geometry: 

\begin{thm}\label{thm:GaussBonnet}
Let $O$ be a Riemannian orbisurface.  Then
\begin{displaymath}
\int_O KdA = 2 \pi \chi (O),
\end{displaymath}
where $K$ is the curvature and $\chi(O)$ is the orbifold Euler
characteristic of $O$.
\end{thm}

\noindent Note that we define the curvature of an orbifold $O$ at a
point $x \in O$ with coordinate chart $(U, \Gamma \backslash
\tilde{U}, \pi)$ to
be the curvature at a lift $\tilde{x} \in \tilde{U}$ of $x$.

Combining the Gauss-Bonnet theorem with Weyl's asymptotic formula, we
see that for an orbisurface with given curvature,
the spectrum determines the orbifold Euler characteristic.  However, since the
orbifold Euler characteristic involves both the genus of the
underlying surface and the orders of the cone points in the
orbisurface, it is not immediately clear that the spectrum determines
the genus.  This is still an open question.

In the case of orientable orbisurfaces, Gordon et al.
\cite{GGWZ} have shown that the Euler characteristic can be recovered from the
asymptotic expansion of the trace of the heat kernel.  Together with
some computations for cone points, this allows them to define a
spectral invariant which determines whether an orbifold is a football
(underlying space $S^2$, cone points at the north and south poles)
or teardrop and determines the orders of the cone points.  In a
similar vein, we give the following obstructions to isospectrality.  

\begin{prop}\label{prop:oneconept}
Fix $g \geq 1$ and $m \geq 2$.  Let $O$ be a compact orientable
Riemann orbisurface of genus $g$ with exactly one cone point of order $m$.
Let $O'$ be in the class of compact orientable Riemann orbisurfaces 
of genus $g$, and suppose that $O$ is isospectral to
$O'$.  Then $O'$ must have exactly one cone point,
and its order is also $m$.   
\end{prop}

\begin{proof}
By Theorem \ref{thm:commoncover}, $O'$
must contain at least one cone point.  We have $\chi (X_{O}) = \chi
(X_{O'})$ by hypothesis, and the observation following Theorem
\ref{thm:GaussBonnet} implies that $\chi(O) = \chi(O')$.  

Suppose that $O'$ has one cone point of order $n_1$.  It follows that 
\begin{displaymath}
\frac{1}{m} = \frac{1}{n_1},
\end{displaymath}
or $m = n_1$.  
Now suppose that $O'$ has two cone points of orders $n_1$ and $n_2$.
Then 
\begin{displaymath}
\frac{1}{m} + 1 = \frac{1}{n_1} + \frac{1}{n_2}.
\end{displaymath}
But $n_i \geq 2$ for $i=1,2$, so $\frac{1}{n_1} + \frac{1}{n_2} \leq
1$.  This is a contradiction, hence $O$ and $O'$ are not isospectral.
This argument is easily extended to the case when $O'$ is assumed to
have more than two cone points.  
\end{proof}
 
We can extend Proposition \ref{prop:oneconept} to the case of two
orbisurfaces with different underlying spaces.

\begin{prop}\label{prop:varyinggenus}
Let $O$ be a compact orientable Riemann orbisurface of genus $g_0 \geq
0$ with $k$ cone points of orders $m_1, \ldots, m_k$, where $m_i \geq
2$ for $i=1,\ldots,k$.  Let $O'$ be a compact orientable Riemann 
orbisurface of genus
$g_1\geq g_0$ with $l$ cone points of orders $n_1, \ldots, n_l$, where
$n_j \geq 2$ for $j=1,\ldots,l$.  Let $h=2(g_0-g_1)$.  If $l \geq
2(k+h)$, then $O$ is not isospectral to $O'$.
\end{prop}

\begin{proof}
Suppose $O$ is isospectral to $O'$.  As in the preceding proof, we
have that $O'$ must contain at least one cone
point.  Also, $\chi (O) = \chi (O')$, i.e.
\begin{displaymath}
2-2g_0-k+\frac{1}{m_1} + \cdots + \frac{1}{m_k} =
2-2g_1-l+\frac{1}{n_1}+\cdots + \frac{1}{n_l}.
\end{displaymath}
So
\begin{displaymath}
\frac{1}{m_1} + \cdots + \frac{1}{m_k} = 2(g_0-g_1)+k-l+\frac{1}{n_1}
+ \cdots + \frac{1}{n_l}
\end{displaymath}
or equivalently
\begin{displaymath}
\frac{1}{m_1} + \cdots + \frac{1}{m_k} \leq h+k-\frac{l}{2}.
\end{displaymath}
But $k+h \leq \frac{l}{2}$ by hypothesis, which gives the desired contradiction.
\end{proof}

\noindent If we assume that $h=0$, i.e. that $g_0=g_1$, then we have
the following special case of Proposition \ref{prop:varyinggenus}.

\begin{corollary}
Fix $g \geq 0$.  Let $O$ be a compact orientable Riemann orbisurface
of genus $g$ with $k$ cone points of orders
$m_1, \ldots, m_k$, $ m_i \geq 2$ for $i=1,\ldots,k$.
Let $O'$ be a compact orientable Riemann orbisurface of genus $g$ with
$l \geq 2k$ cone
points of orders $n_1, \ldots, n_l$, $n_j \geq 2$ for $j=1, \ldots,
l$.  Then $O$ is not isospectral to $O'$.
\end{corollary}

Note that in all of the above results, we have the hypothesis that
$O'$ is hyperbolic.  In the usual case of surfaces, we know that any
surface isospectral to a given one of fixed constant curvature must
have the same constant curvature.  The proof of this uses the
asymptotic expansion of the trace of the heat kernel; in the case of
orbisurfaces, this expansion is more complicated and it is no longer
clear that fixed constant curvature is a spectral invariant.

\section{Huber's Theorem}\label{sec:hubers}

Huber's theorem is a powerful tool in the study of questions of
isospectrality of compact Riemann surfaces.  It allows us to translate
information about eigenvalues into information about the geometry of
the surface, and specifically about the lengths of closed geodesics on
the surface (see \cite{Bubook}).

\begin{thm}\label{thm:hubers}
(Huber) Two compact Riemann surfaces of genus $g \geq 2$ have the same
spectrum of the Laplacian if and only if they have the same length spectrum.
\end{thm}

\noindent The length spectrum is the sequence of all
lengths of all oriented closed geodesics on the surface, arranged in
ascending order.

The idea of the proof is as follows.  First, a fundamental domain
argument leads to a length trace formula.  The known eigenfunction
expansion of the heat kernel as given in (\ref{eqn:heatkernel}) is then
plugged into this length trace formula to obtain the Selberg Trace
Formula.  The Selberg Trace Formula contains information about the
eigenvalues on one side and information about the lengths of closed
geodesics and the area of the Riemann surface on the other; it is then
a matter of showing that the eigenvalues determine the lengths of the closed
geodesics and vice versa.   Further background and the case of the
Selberg Trace Formula for compact hyperbolic manifolds can be found in 
Randol (in \cite{Chavelbook}).  Other sources for the Riemann surface
case include \cite{Bubook}, \cite{Hejhalbook1} and \cite{Pitkin}.

We want to extend Huber's theorem to the class of compact orientable Riemann
orbisurfaces.  To begin, we need to
exhibit a Selberg Trace Formula for such objects.  
We can state this formula in terms of a function
$h(r)$ and its Fourier transform $g(u)$, which are not required to have compact support;
namely, $h(r)$ satisfies the following (weaker) conditions:
\begin{assumption}\label{assump:Hejhal7.2} 
\
\begin{itemize}
\item $h(r)$ is an analytic function on $|Im (r)| \leq \frac{1}{2} + \delta $;
\item $h(-r) = h(r)$;
\item $|h(r)| \leq M[1+ |Re(r)|]^{-2 -\delta}$.
\end{itemize}
\noindent The numbers $\delta$ and $M$ are some positive constants.
\end{assumption}

Hejhal \cite{Hejhalbook1} obtains the following version of the Selberg Trace
Formula for the case of interest:
\begin{thm}\label{thm:ofldSTF} 
Suppose that 
\begin{itemize}
\item $\Gamma \subset PSL(2, \mathbb{R})$ is a Fuchsian group with
  compact fundamental region;
\item $h(r)$ satisfies Assumption \ref{assump:Hejhal7.2};
\item $\{ \phi_n \}_{n=0}^{\infty}$ is an orthonormal eigenfunction
  basis for $L^2 ( \Gamma \backslash \mathbb{H})$.
\end{itemize}
Then
\begin{eqnarray}\label{eqn:ofldSTF}
\sum_{n=0}^{\infty} h(r_n) & = & \frac{\mu(F)}{4 \pi}
\int_{-\infty}^{\infty} r h(r) \tanh (\pi r) dr \nonumber \\
& & \mbox{}
+ \sum_{\genfrac{}{}{0pt}{}{\{ R \}}{\text{elliptic}}} \frac{1}{2 m(R)
  \sin \theta (R) } \int_{-\infty}^{\infty} \frac{e^{-2 \theta (R)
    r}}{1 + e^{-2 \pi r}} h(r) dr \nonumber \\
& & \mbox{} + \sum_{\genfrac{}{}{0pt}{}{\{ P \}}{\text{hyperbolic}}}
\frac{\ln N(P_c)}{N(P)^{1/2} - N(P)^{-1/2}} g[\ln N(P)], 
\end{eqnarray}
where all the sums and integrals in sight are absolutely convergent.
\end{thm}

Note that the left side of (\ref{eqn:ofldSTF}) encodes information about the 
eigenvalues of the Laplacian, as $\lambda_n = -\frac{1}{4} - r_n^2$. 
(For easier reference, we follow Hejhal's convention of nonpositive eigenvalues 
for this section.) 
The first and last terms on the right side of (\ref{eqn:ofldSTF}) match the
terms which appear in the Selberg Trace Formula for compact Riemann surfaces,
where $\mu(F)$ is the area of a fundamental domain for $\Gamma$ and 
the last term is the sum over hyperbolic conjugacy classes.
Every hyperbolic element $P$ in $PSL(2, \mathbb{R})$ is conjugate to a
dilation $z \mapsto m^2 z, m>1$; we call $m^2$ the \emph{norm} of $P$ and
denote it by $N(P)$.  A \emph{primitive} hyperbolic element is one which cannot 
be written as a nontrivial power of another hyperbolic element.  In (\ref{eqn:ofldSTF}), 
$P_c$ is a primitive hyperbolic element with $P= P_c^k$ for some $k \geq 1$.  
Note that \cite[Prop. 2.3]{Hejhalbook1}
$$
{\genfrac{}{}{0pt}{}{\text{inf}}{z \in \mathbb{H}}} \ \ d(z, Tz) = \ln N(T),
$$
where $d$ is the distance on $\mathbb{H}$.  The infimum is realized by
all points $z$ which lie on the geodesic in $\mathbb{H}$ which is
invariant under the action of $T$.
In particular, we allow geodesics on $O = \Gamma \backslash \mathbb{H}$ 
to pass through cone points, and include the lengths of such geodesics 
in the length spectrum of $O$.  In the sum over elliptic conjugacy
classes, $m(R)$ denotes the order of the centralizer (in $\Gamma$) of a
representative $R$ and $\theta (R)$ is half the angle of rotation.
We have $Tr(R) = 2 \cos \theta(R)$, and $0 < \theta < \pi$.   

In the proof of our partial extension of Huber's theorem, we will need
the following result of Stanhope \cite{Stanhopebounds}.

\begin{thm}\label{thm:Stanhopefiniteisotropy}
Let $\mathcal{S}$ be a collection of isospectral orientable compact
Riemannian orbifolds that share a uniform lower bound $\kappa (n-1)$,
$\kappa$ real, on Ricci curvature.  Then there are only finitely many
possible isotropy types, up to isomorphism, for points in an orbifold
in $\mathcal{S}$.
\end{thm}

We are now prepared to state our partial extension of Huber's theorem
to the setting of compact orientable Riemann orbisurfaces.

\begin{thm}\label{thm:partialhubers}
If two compact orientable Riemann orbisurfaces are Laplace
isospectral, then we can determine their
length spectra and the orders of their cone points,
up to finitely many possibilities.  Knowledge of the
length spectrum and the orders of the cone points determines the
Laplace spectrum.
\end{thm}

\begin{proof}  
We will consider Theorem \ref{thm:ofldSTF} for a specific function
$h(r)$. 
Fix $t>0$ and let $h(r) = e^{-r^2 t}$.  Then $h(r)$ satisfies
Assumption \ref{assump:Hejhal7.2},  and we have
\begin{eqnarray*}
g(u) & = & \frac{1}{2 \pi} \int_{- \infty}^{\infty} h(r) e^{-iru} dr \\
& = & \frac{1}{\sqrt{4 \pi t}}e^{-u^2/4t}
\end{eqnarray*}  
where the first line is the definition of $g(u)$ as the Fourier
transform of $h(r)$ and the second line follows from Fourier analysis
using a standard polar coordinates trick.

By Theorem \ref{thm:ofldSTF}, we have
\begin{eqnarray}\label{eqn:ofldSTFh(r)}
\sum_{n=0}^{\infty} e^{-r_n^2 t} & = & \frac{\mu(F)}{4 \pi}
\int_{-\infty}^{\infty} r e^{-r^2 t} \tanh (\pi r) dr \nonumber \\
& & \mbox{}
+ \sum_{\genfrac{}{}{0pt}{}{\{ R \}}{\text{elliptic}}} \frac{1}{2 m(R)
  \sin \theta (R) } \int_{-\infty}^{\infty} \frac{e^{-2 \theta (R)
    r}}{1 + e^{-2 \pi r}} e^{-r^2 t} dr \nonumber \\
& & \mbox{} + \sum_{\genfrac{}{}{0pt}{}{\{ P \}}{\text{hyperbolic}}}
\frac{\ln N(P_c)}{N(P)^{1/2} - N(P)^{-1/2}}  \frac{1}{\sqrt{4 \pi
    t}}e^{-(\ln N(P))^2/4t}. 
\end{eqnarray}

Let $O$ and $O'$ be compact orientable Riemann orbisurfaces with the
same Laplace spectrum. 
By Theorem \ref{thm:Weylasymp}, the Laplace spectrum determines an
orbifold's volume.  So we must have vol($O$) = vol($O'$), and thus the
first term on the right side of (\ref{eqn:ofldSTFh(r)}) must be the same for
$O$ and $O'$.  
  
Note that both $O$ and $O'$ have a metric of constant curvature -1, 
and hence share a uniform lower bound on their Ricci curvature.
By Theorem \ref{thm:Stanhopefiniteisotropy}, we know that there are only
finitely many possible isotropy types, up to isomorphism, for points
in $O$ [$O'$].
It is well-known (e.g. \cite[p. 16]{FischerSTF}) that there are
only finitely many cone points in $O$ [$O'$].
Putting these two facts together, we see that the Laplace spectrum 
determines (up to finitely many possibilities) the orders of the cone 
points in $O$ [$O'$], and thus the sum
over the elliptic elements in (\ref{eqn:ofldSTFh(r)}).
That is, up to finitely many possibilities, we
know the function   
\begin{displaymath}
f(t) = \sum_{\genfrac{}{}{0pt}{}{\{ P \}}{\text{hyperbolic}}}
\frac{\ln N(P_c)}{N(P)^{1/2} - N(P)^{-1/2}} e^{-(\ln N(P))^2/4t}.
\end{displaymath}

Consider the function $f(t)e^{\omega^2/4t}$.
Take the limit of this function as $t \downarrow 0$; we see that there
is a unique $\omega >0$ for which this limit is finite and nonzero.
Let $\gamma_1$ be the shortest primitive closed geodesic in $O$.  Then
$\omega = \ell(\gamma_1)$.  We remove the contribution of
$\gamma_1$ and all its powers from $f(t)$, and proceed as above to
find the length of the next-shortest primitive closed geodesic.  In
this way, we can determine the lengths of the hyperbolic elements in
$O$ [$O'$], up to finitely many possible lists of lengths.   

Now suppose we know the length spectrum and the orders of the cone
points for $O$ [$O'$].  The argument that we then know the Laplace
spectrum of $O$ [$O'$] is as for Riemann surfaces (see
\cite[p. 45] {Pitkin}).  We include it for completeness.  

First, we multiply both sides of (\ref{eqn:ofldSTFh(r)}) by $e^{-t/4}$ and
recall that $\lambda_n = -\frac{1}{4} - r_n^2$ to obtain
\begin{eqnarray}\label{eqn:ofldSTFusefulform}
\sum_{n=0}^{\infty} e^{\lambda_n t} & = & \frac{\mu (F)}{4 \pi}
e^{-t/4} \int_{-\infty}^{\infty} r e^{-r^2 t} \tanh (\pi r) dr \nonumber \\
& & \mbox{} + \sum_{\genfrac{}{}{0pt}{}{\{ R \}}{\text{elliptic}}}
\frac{e^{-t/4}}{2 m(R) \sin \theta (R) } \int_{-\infty}^{\infty} \frac{e^{-2
    \theta (R) r}}{1 + e^{-2 \pi r}} e^{-r^2 t} dr \nonumber \\
& & \mbox{} + \sum_{\genfrac{}{}{0pt}{}{\{ P \}}{\text{hyperbolic}}}
\frac{\ln N(P_c)}{N(P)^{1/2} - N(P)^{-1/2}} \cdot
\frac{e^{-t/4}}{\sqrt{4 \pi t}} e^{-(\ln N(P))^2/4t}.
\end{eqnarray} 
Knowledge of
the length spectrum and the orders of the cone points in
(\ref{eqn:ofldSTFusefulform}) translates to knowledge of the function
\begin{eqnarray*}
c(t) & = & \sum_{n=0}^{\infty} e^{\lambda_n t} - \frac{\mu (F)}{4 \pi}
e^{-t/4} \int_{-\infty}^{\infty} r e^{-r^2 t} \tanh (\pi r) dr \nonumber \\
& = & \sum_{-\frac{1}{4} \leq \lambda_n < 0} e^{\lambda_n t} -
\frac{\mu (F)}{4 \pi} e^{-t/4} \int_{-\infty}^{\infty} r e^{-r^2 t}
\tanh (\pi r) dr + \sum_{\lambda_n < -\frac{1}{4}} e^{\lambda_n t} \nonumber \\
& = &  \sum_{-\frac{1}{4} \leq \lambda_n < 0} e^{\lambda_n t} - \sigma
(t) e^{-t/4} \mu (F) + \sum_{\lambda_n < -\frac{1}{4}} e^{\lambda_n t},
\end{eqnarray*}
where
\begin{displaymath}
\sigma (t) = \frac{1}{2 \pi} \int_0^{\infty} r e^{-r^2 t} \tanh (\pi r) dr.
\end{displaymath}
It is a straightforward but tedious calculation to show that as $t \rightarrow \infty$, $\sigma (t) \rightarrow 0$ (see \cite{EBDthesis}).  

If $\lambda_1 \geq -\frac{1}{4}$, then $-\lambda_1$ is the unique
$\omega > 0$ such that 
\begin{displaymath}
0 < \lim_{t \rightarrow \infty} e^{\omega t} c(t) < \infty.
\end{displaymath}  
In fact, this limit is the multiplicity $m_1$ of $\lambda_1$.  We can
therefore subtract $m_1 e^{\lambda_1 t}$ from $c(t)$ and continue in
this way to find all the small eigenvalues.  Once all the small
eigenvalues have been found, the function
\begin{displaymath}
\tilde{c} (t) = - \sigma (t) e^{-t/4} \mu (F) + \sum_{\lambda_n <
  -\frac{1}{4}} e^{\lambda_n t}
\end{displaymath}
has the property that for $\omega >0$,
\begin{displaymath}
\lim_{t \rightarrow \infty} e^{\omega t} \tilde{c} (t)
\end{displaymath}
is $0$ or $\infty$.  So we can now multiply $\tilde{c} (t)$ by $-
\frac{e^{t/4}}{\sigma (t)}$ and take the limit as $t \rightarrow \infty$
to get $\mu (F)$.  We then know the function
\begin{displaymath}
\sum_{\lambda_n < -\frac{1}{4}} e^{\lambda_n t},
\end{displaymath}
and we can determine the remaining eigenvalues in the same way as we
found the small eigenvalues.
Hence the spectrum of the Laplacian is determined by the length
spectrum and the orders of the cone points.
\end{proof}

\section{Finiteness of Isospectral Sets}\label{sec:finiteness}

McKean \cite{McK72} showed that only finitely many compact Riemann
surfaces have a given spectrum.  We extend this result to the setting
of compact orientable Riemann orbisurfaces.  Specifically, we show

\begin{thm}\label{thm:finitelymany}
Let $O$ be a compact orientable Riemann orbisurface of genus $g
\geq 1$.  In the class
of compact orientable hyperbolic orbifolds, there are only finitely many members which are
isospectral to $O$.
\end{thm}

\begin{remark}
Note that there is no need for a dimension restriction on the
orbifolds that can be isospectral to $O$, by Theorem
\ref{thm:Weylasymp}.  In addition, by Theorem \ref{thm:commoncover},
there can be no Riemann surfaces isospectral to $O$.
\end{remark}

McKean \cite{McK72} states the following proposition, which he attributes to
Fricke and Klein \cite{FrKlein}.  For the sake of accuracy, we give a complete proof below.    

\begin{prop}\label{prop:tracesdetermine}
Let $M = G \backslash \mathbb{H}$ be a Riemann surface of genus $g \geq 2$, 
where $G \leq SL(2,\mathbb{R})$. Let the set \{$h_1,\ldots, h_n$\} 
, $n \leq 2g$, be a generating set for $G$.  Then the single, double,
and triple traces 
\begin{eqnarray*}
tr(h_i) ,& \\
tr(h_i h_j) ,& i<j \\
tr (h_i h_j h_k) ,& i<j<k
\end{eqnarray*} 
determine $G$ up to a conjugation in $PSL(2,\mathbb{R})$ or a reflection.
\end{prop}

\begin{proof}
Let $G$ and $G'$ be two subgroups of $SL(2,\mathbb{R})$ with the same
single, double, and triple traces of their generators.  
Fix $h_1 \in G$.  
Since the single traces of the generators of $G$ and $G'$ are equal, 
we can pair $h_1$ with an element in $G'$ that translates the same amount; 
that is, we can suppose that $h_1 = h_1'$ and
that $h_1(z) = m^2 z$ with $m>1$.  
Note that any other diagonal element in $G$ fixes the same geodesic in
$\mathbb{H}$ as $h_1$ and is thus a multiple of $h_1$.  So we can
assume that $h_1$ is the only diagonal element in $\{h_1, \ldots, h_n\}$.
For $i>1$, we have 
\begin{equation}\label{eqn:traceh_i}
tr(h_i) = a_i + d_i = tr(h_i') = a_i' + d_i'.
\end{equation}
Also,
\begin{displaymath}
tr(h_1h_i) = tr \left( 
\begin{pmatrix}
m & 0 \\ 
0 & m^{-1}  
\end{pmatrix}
\begin{pmatrix}
a_i & b_i \\
c_i & d_i
\end{pmatrix}
\right) = tr \left(
\begin{pmatrix}
ma_i & mb_i \\
m^{-1}c_i & m^{-1}d_i
\end{pmatrix}
\right) = 
ma_i + m^{-1}d_i
\end{displaymath}
and similarly for $tr(h_1'h_i')$, so that
\begin{displaymath}
ma_i + m^{-1}d_i = tr(h_1h_i) = tr(h_1'h_i') = ma_i' + m^{-1}d_i',
\end{displaymath}
or equivalently
\begin{equation}\label{eqn:mdiffs}
m(a_i - a_i') = m^{-1}(d_i' - d_i).
\end{equation}
From equation (\ref{eqn:traceh_i}) we see that $a_i - a_i' =
d_i'-d_i$.  But we assumed $m>1$, so we must have 
\begin{equation}\label{eqn:diagentriesequal}
a_i = a_i' \text{   and   } d_i = d_i'.
\end{equation}  
We also know that det($h_i$) = det ($h_i'$) = 1 for
all $i$, so 
\begin{equation}\label{eqn:offdiags}
b_ic_i = b_i'c_i' 
\end{equation}
for all $i$.  Straightforward calculations
show that
\begin{equation}\label{eqn:doubletrdiff}
tr(h_ih_j)-tr(h_i'h_j') = b_ic_j-b_i'c_j'+c_ib_j - c_i'b_j'
\end{equation}
and
\begin{equation}\label{eqn:tripletrdiff}
tr(h_1h_ih_j) - tr(h_1'h_i'h_j') = m(b_ic_j - b_i'c_j')+m^{-1}(c_ib_j
- c_i'b_j')
\end{equation}
for $1<i<j$.  Combining equations (\ref{eqn:doubletrdiff}) and
(\ref{eqn:tripletrdiff}) as we combined (\ref{eqn:traceh_i}) and 
(\ref{eqn:mdiffs}), we see that
$b_ic_j = b_i'c_j'$ for $1<i<j$.  We want to see that none of these
numbers are zero.  Suppose $c_2=0$.  Then 
\begin{displaymath}
h_1^{-n}h_2h_1^n(\sqrt{-1}) = \frac{a_2 \sqrt{-1} +m^{-2n}b_2}{d_2},
\end{displaymath}
where this is the M\"{o}bius action of $SL(2, \mathbb{R}$) on $\mathbb{H}$.
So we get infinitely many images of $\sqrt{-1}$ accumulating at
$\frac{a_2}{d_2}\sqrt{-1} \in \mathbb{H}$ (unless $b_2=0$, which
implies that $h_2$ is diagonal, contradicting our assumption that
$h_1$ is the only diagonal element in $\{h_1, \ldots, h_n\}$).  
This contradicts the assumption
that $G$ acts properly discontinuously on $\mathbb{H}$.  
A similar argument with $b_2 = 0$ shows that the off-diagonal entries in 
the matrix representing the element $h_2$ are nonzero.  
Our choice of $h_2$ was arbitrary,  
thus none of the off-diagonal entries in the matrices representing the
elements $h_2,\ldots,h_n$ and $h_2',\ldots,h_n'$ are zero.  We have
\begin{displaymath}
\frac{c_j'}{c_j} = \frac{b_i}{b_i'} = \frac{c_i'}{c_i},
\end{displaymath}
where the second equality is equation (\ref{eqn:offdiags}), and this
common ratio is independent of $i>1$.  Since the traces do
not tell us whether the ratio is positive, we must allow the
reflection
\begin{displaymath}
G \rightarrow 
\begin{pmatrix}
1 & 0 \\
0 & -1
\end{pmatrix}
G
\begin{pmatrix}
1 & 0 \\
0 & -1
\end{pmatrix}.
\end{displaymath}

Suppose that the common ratio is equal to $t^2$, i.e.
\begin{displaymath}
\frac{b_i}{b_i'} = t^2 = \frac{c_i'}{c_i}
\end{displaymath}
for all $i>1$.  Then $b_i = t^2 b_i'$ and $c_i = t^{-2}c_i'$ for $i>1$.  Thus
\begin{displaymath}
\begin{pmatrix}
a_i & b_i \\
c_i & d_i
\end{pmatrix}
=
\begin{pmatrix}
a_i & t^2 b_i' \\
t^{-2}c_i' & d_i
\end{pmatrix}
\end{displaymath}
for all $i>1$.  We saw that $a_i = a_i'$ and $d_i = d_i'$ for all $i>1$, 
so
\begin{displaymath}
h_i = 
\begin{pmatrix}
a_i & b_i \\
c_i & d_i
\end{pmatrix}
=
\begin{pmatrix}
a_i' & t^2 b_i' \\
t^{-2}c_i' & d_i'
\end{pmatrix}
=
s h_i' s^{-1}
\end{displaymath}
for $s \in SL(2, \mathbb{R})$ given by $s: z \mapsto t^2 z$.  
Thus there exists $s \in SL(2, \mathbb{R})$ which, for all $i$,
conjugates $h_i$ to $h_i'$.
Hence $G$ and $G'$
are the same group up to conjugation in $PSL(2, \mathbb{R})$ or a reflection. 
\end{proof}

Note that we can easily extend this result to the case of a group $G$
which is the fundamental group of a compact orientable Riemann orbisurface of genus
$g \geq 1$.  We know that
any such group contains a hyperbolic element; without loss of
generality, label this element $h_1$.  Then the calculations which
show that equation (\ref{eqn:diagentriesequal}) holds are still valid,
as are the calculations which show that $b_ic_j = b_i'c_j'$ for $i
\neq j$.  An elliptic element $R$ in
$PSL(2, \mathbb{R})$ is conjugate to an element of the form
\begin{displaymath}
\begin{pmatrix}
\cos \theta & -\sin \theta \\
\sin \theta & \cos \theta
\end{pmatrix}
\end{displaymath} 
for $0 < \theta < \pi$; a straightforward but messy calculation shows 
that the off-diagonal entries of $R$ are zero only if $\theta = 0, \pi$.  
So the remainder of the argument holds in the desired setting.

To prove Theorem \ref{thm:finitelymany}, we will need the following
result of Stanhope \cite{Stanhopebounds} which gives an upper bound
on the diameter of an orbifold:

\begin{prop}\label{prop:Stanhopediameter}
Let $O$ be a compact Riemannian orbifold with Ricci curvature bounded
below by $\kappa (n-1), \kappa$ real.  Fix an arbitrary constant $r$
greater than zero.  Then the number of disjoint balls of radius $r$
that can be placed in $O$ is bounded above by a number that depends
only on $\kappa$ and the number of eigenvalues of $O$ less than or
equal to $\lambda_{\kappa}^n(r)$.  

In particular the diameter of $O$ is bounded above by a number that
depends only on Spec($O$) and $\kappa$.
\end{prop}

\noindent So isospectral families of orbifolds with a common uniform lower Ricci
curvature bound also have a common upper diameter bound. 

In Beardon \cite{Beardonbook}, a careful study of the geometry of discrete
groups is undertaken.  In particular, he looks at Fuchsian groups, which may be
considered as discrete groups of isometries of the hyperbolic plane.
To every Fuchsian group $G$, we can associate its
Dirichlet polygon; a Dirichlet polygon centered at a point $w \in
\mathbb{H}$ is the set of all points which are, among all their images
under the action of $G$, closest to $w$.  Beardon proves the following
theorem about such a polygon $P$ \cite[Thm. 9.3.3]{Beardonbook}:

\begin{thm}\label{thm:Beardonsidepairing}
The set of side-pairing elements $G^*$ of $P$ generate $G$.
\end{thm}

\noindent Thus every Fuchsian group $G$ is generated by the side-pairing
elements of a Dirichlet polygon associated to it.  
Note that we can regard the elliptic fixed points of
$G$ as vertices of $P$.  

We are now ready to prove Theorem \ref{thm:finitelymany}.  

\begin{proof}[Proof of Theorem \ref{thm:finitelymany}]
Let $C$ be the class of compact orientable hyperbolic orbifolds, and let $S$ denote 
the subclass of $C$ containing those
orbifolds which are isospectral to $O$.  Note that any member of $S$
is a compact orientable Riemann orbisurface, and that it is determined by
its fundamental group $\Gamma$.  By Proposition \ref{prop:tracesdetermine}, 
specifying $\Gamma$ (up to a reflection or
conjugation) is the same as specifying the single, double and
triple traces of a set of generators.  
We want to show that the spectrum determines finitely many possibilities for 
these traces, thus only finitely many choices for $\Gamma$,
and hence that $S$ is a finite set.

Theorem \ref{thm:partialhubers} tells us that the Laplace spectrum
determines (up to finitely many possibilities) the length $\ell(Q)$ of
a shortest closed path in the free homotopy class associated to a
given hyperbolic conjugacy class $Q$ in $\Gamma$.  The relation
between the trace of $Q$ and $\ell(Q)$ is given by:
\begin{displaymath}
tr(Q) = \pm 2 \cosh \frac{1}{2} \ell(Q).
\end{displaymath}   
Note that $\frac{1}{2} \ell(Q)$ is bounded by $D$, the diameter of 
$\Gamma \backslash \mathbb{H}$.
Thus the single traces of the hyperbolic conjugacy classes are
bounded by $2\cosh D$.  
By choosing our generating set to be the set of side-pairing elements 
of a Dirichlet polygon $P$ for $\Gamma$, we can easily bound the traces 
of hyperbolic conjugacy classes which arise as a product of two or three generators.
We fix a point $p \in P$ and determine an
upper bound for $\text{dist} (p, g_2 \circ g_1 (p))$, where $g_1$
and $g_2$ are side-pairing elements of $P$ and $g_2
\circ g_1$ is hyperbolic.  We have
$$   
\text{dist}(p, g_2 \circ g_1 (p)) \leq \text{dist}(p, g_1 (p)) +
\text{dist}(g_1(p), g_2 \circ g_1 (p))
$$
and each term on the right side is bounded by $2D$.  Thus 
$$
tr(g_2 \circ g_1) \leq 2 \cosh 2D.
$$
A similar argument shows that the trace of the product of three
side-pairing elements which is hyperbolic is bounded by $2 \cosh 3D$.

By Proposition \ref{prop:Stanhopediameter}, there is
a common upper bound on the diameter of any orbisurface in $S$.  
So there are
only finitely many possibilities for the 
trace of a hyperbolic element which arises as a product of one, two or three
side-pairing elements of $P$.    

Finally, we need to consider the case in which a side-pairing element of 
$P$ (or a product of such elements) is elliptic.  
Beardon \cite[p.225]{Beardonbook} notes that $P$
contains representatives of all conjugacy classes of elliptic elements
in $\Gamma$.    
But we know that $\Gamma \backslash \mathbb{H}$ contains only finitely many 
cone points, and Theorem \ref{thm:Stanhopefiniteisotropy}
tells us that there can be only finitely many possible isotropy types for
the points in $\Gamma \backslash \mathbb{H}$.   So there are only
finitely many choices for the trace of any elliptic element in
$\Gamma$; this implies that if a product of generators of $\Gamma$ is elliptic,
there are only finitely many choices for the trace of such a product. 

Thus we have shown that there are only finitely many ways to choose
the generators of $\Gamma$, up to a reflection or conjugation. 
\end{proof}

\bibliographystyle{plain}
\bibliography{../emily}

\end{document}